\documentclass{amsart}
\usepackage{amsthm}
\usepackage{enumerate}

\makeatletter
\newtheorem*{rep@theorem}{\rep@title}
\newcommand{\newreptheorem}[2]{%
\newenvironment{rep#1}[1]{%
 \def\rep@title{#2 \ref{##1}}%
 \begin{rep@theorem}}%
 {\end{rep@theorem}}}
\makeatother

\newtheorem{theorem}{Theorem}
\newreptheorem{theorem}{Theorem}
\newtheorem{lemma}[theorem]{Lemma}

\newtheorem{conjecture}{Conjecture}
\newtheorem{question}{Question}
\newtheorem{proposition}[theorem]{Proposition}

\theoremstyle{definition}
\newtheorem{definition}{Definition}
\newtheorem{remark}{Remark}

\newtheoremstyle{named}{}{}{\itshape}{}{\bfseries}{.}{.5em}{\thmnote{#3}}
\theoremstyle{named}
\newtheorem*{namedtheorem}{}

\newcommand{\II}{\operatorname{II}}
\newcommand{\dHaus}{\operatorname{d\mathcal{H}}}
\newcommand{\osc}{\operatorname*{osc}}

\newcommand{\dvol}{\operatorname{dvol}}

\title{Singularities of the Lagrangian Mean Curvature Flow}
\date{}

\author{Andrew A Cooper}

\address{Box 8205, North Carolina State University, Raleigh, NC 27695}
	
\email{andrew.cooper@math.ncsu.edu}
\thanks{The author wishes to thank Jon Wolfson, Jerry Kazdan, and Christopher Croke for their support and encouragement of this project, and Radmila Sazdanovi\'c for her support and comments. He also wishes to thank the Park City Mathematics Institute, where parts of this manuscript were prepared. Theorem A appeared in the author's Ph.D.~thesis \cite{cooperthesis}.}

\begin{document}
\maketitle

\begin{abstract}
In this paper we investigate the singularities of Lagrangian mean curvature flows in $\mathbf{C}^m$ by means of smooth singularity models. Type I singularities can only occur at certain times determined by invariants in the cohomology of the initial data. In the type II case, these smooth singularity models are asymptotic to special Lagrangian cones; hence all type II singularities are modeled by unions of special Lagrangian cones.
\end{abstract}

\section{Introduction}

Mean curvature flow, the downward gradient flow for the area functional on submanifolds, has been extensively studied in codimension one. Mean curvature flow of Lagrangian submanifolds in K\"ahler-Einstein manifolds is among simplest settings for high-codimension mean curvature flow.  In case the ambient manifold is $\mathbf{C}^m$, compact flows must become singular in finite time.  Conjectures of Thomas-Yau and Wang (see the article by Joyce \cite{joycesurvey} for a thorough discussion) relate Lagrangian mean curvature flow to the structure of Hamiltonian isotopy classes of Calabi-Yau manifolds. Neves \cite{nevessurvey} and Wolfson \cite{wolfson05} have shown that Lagrangian mean curvature flow singularities can be expected to interfere with this approach, even generically. Thus singularity analysis for Lagrangian mean curvature flow is required.

In this paper we prove that mild, scale-respecting singularities (so-called type I) can only occur at certain times depending only on the Lagrangian invariants of the initial data, and that all other singularities are modeled by (collections of) special Lagrangian cones. The latter result supports Joyce's conjectured resolution of the Thomas-Yau-Wang conjectures by Lagrangian mean curvature flow with conic singularities, and may also enable Lagrangian mean curvature flow with surgery.  Progress toward flow with isolated conic singularities has been made by Behrndt \cite{behrndt11}, under the assumption that the cones are stable special Lagrangian. Begley-Moore \cite{begley&moore15} have shown that Lagrangians with a finite number of cone singularities, each modeled on a pair of transversally intersecting Lagrangian planes, can be taken as initial data for a smooth Lagrangian flow.

\subsection{Summary of Results}
The central observation of this paper is:
	\begin{namedtheorem}[Theorem A]
		Let $\Sigma(t)$ be a Lagrangian mean curvature flow defined on a maximal interval $[0,T)$, $T<\infty$. There exists a smooth model of the singularity at time $T$ which is:
			\begin{enumerate}
				\item in the type I case, a monotone Lagrangian.
				\item in the type II case, an exact, zero-Maslov Lagrangian.
			\end{enumerate}
		Moreover the 1-cycles which survive to this smooth singularity model are predicted by the Liouville and Maslov classes of the initial data $\Sigma(0)$.
	\end{namedtheorem}

From Theorem A, we get the following two  important corollaries regarding type II singularities of the Lagrangian mean curvature flow.
	\begin{namedtheorem}[Theorem B]
		Any type II singularity of Lagrangian mean curvature flow is modeled by a union of special Lagrangian cones.
	\end{namedtheorem}
	
	\begin{namedtheorem}[Theorem B']
		Any type II singularity of mean curvature flow starting from an almost-calibrated zero-Maslov Lagrangian is modeled by a smooth special Lagrangian.
		
		Any type II singularity of mean curvature flow starting from a monotone Lagrangian surface is modeled by a smooth special Lagrangian.
	\end{namedtheorem}

	We can also use our smooth singularity models to understand something about the type I case. 
	\begin{namedtheorem}[Theorem C]
		Let $\Sigma$ be a compact Lagrangian submanifold of $\mathbf{C}^m$. If mean curvature flow starting from $\Sigma$ has singular time $T$, and the singularity is of type I,  then $T$ has the form $\frac{\lambda.\gamma}{2h.\gamma}$, where $[\lambda]$ and $[h]$ are the Liouville and Maslov classes of $\Sigma$, and $\gamma$ is some cycle in $\Sigma$.
	\end{namedtheorem}

\subsection{Previous work}

	Previously, Chen-Li, Neves and Groh-Schwarz-Smoczyk-Zehmisch investigated the nature of singularities of Lagrangian mean curvature flows in $\mathbf{C}^m$ in case the initial data are zero-Maslov or monotone  \cite{groh&schwarz&smoczyk&zehmisch} \cite{neveszeromaslov} \cite{nevesmonotone} . As we will see in Theorem A, these two cases essentially capture the singular behavior of general Lagrangian mean curvature flows.
	
	\begin{theorem}[\cite{neveszeromaslov}]\label{TF.zeromaslov}
		The tangent flow to a zero-Maslov Lagrangian mean curvature flow is a collection of minimal cones. Moreover, the convergence of the tangent flow rescaling sequence respects the Lagrangian angle: each minimal cone $\mu_k$ has an associated $\beta_k\in\mathbf{R}$ so that for any $\phi\in C_0^\infty(\mathbf{C}^m)$ and any $f\in C^\infty(\mathbf{R})$, 
			\begin{equation}
				\sum_{k=1}^N m_k\mu_k(\phi)f(\beta_k)=\lim_j\int_{\tilde{L}_j}\phi f(\beta)\dHaus^m
			\end{equation}
			
		In case the initial data are almost-calibrated, $N=1$.
	\end{theorem}
	
	\begin{theorem}[\cite{groh&schwarz&smoczyk&zehmisch},\cite{nevesmonotone}]\label{TF.monotone}
		Suppose $L:\Sigma\times [0,T)\rightarrow \mathbf{C}^m$ is a compact Lagrangian mean curvature flow whose initial data are monotone, i.e. $[\lambda(0)]=C[h(0)]$. If $T<\frac{1}{2}C$, then the tangent flow to $L$ is a collection of minimal cones. Moreover, the convergence respects the Lagrangian angle.
		
		In the case $m=2$, $N=1$.
	\end{theorem}
	
	Theorem \ref{SBU.SS} shows that in the cases above, the singularity must have been of type II.

\subsection{Plan of the Paper}
	In Section 2  we introduce the basics of the Lagrangian mean curvature flow, including the evolution of Lagrangian invariants and singularity analysis via rescaling.
	
	In Section 3 we prove Theorem A. In Section 4 we prove a ``maximum principle" for $S^1$-valued functions satisfying the heat equation. In Section 5, we apply this maximum principle to show that singularities which are modeled by zero-Maslov smooth blow-ups admit minimal-cone singularity models, which implies theorem B. In Sections 6, we prove Theorem C.   In Section 7 we consider the situations in Theorems \ref{TF.zeromaslov} and \ref{TF.monotone} in which the tangent flow has a single Lagrangian angle and prove Theorem B'.  In section 8, we outline some questions for future work independent of the Thomas-Yau conjecture.

  \section{Preliminaries}

\subsection{Lagrangian geometry along the flow}

\begin{definition}
	A {\em Liouville form} for the K\"{a}hler manifold $(X,\omega,J)$ is any one-form $\eta$ on $X$ with $d\eta=2\omega$. We call $\lambda=L^*\eta$ the {\em Liouville form} of $L$.
\end{definition}

Note that if $L:\Sigma\rightarrow X$ is a Lagrangian submanifold, then $\lambda=L^*\eta$ is closed, since $L^*\omega=0$.

\begin{definition}
	The Maslov form $h$ is the one-form dual, with respect to $\omega$, to the mean curvature vector $H=H^i\nu_i$, that is, $h=L^*(H\lrcorner\omega)$.
\end{definition}

When $(X,\omega,J)$ is K\"{a}hler-Einstein, by Codazzi's equation and the contracted Bianchi identity, we have that $\nabla_iH_j=\nabla_jH_i$, hence $h$ is closed. 

The following lemma allows us to consider {\em Lagrangian}  mean curvature flow.

\begin{lemma}
	Mean curvature flow preserves the Lagrangian condition in a K\"{a}hler-Einstein manifold, i.e. if the initial submanifold $L_0$ is Lagrangian, so is each time slice $L_t$.
\end{lemma}
\begin{proof}
	The Lagrangian condition is $L^*\omega=0$. We compute
	\begin{equation}
		\begin{aligned}
			\frac{\partial}{\partial t}L^*\omega&=L^*(\mathcal{L}_{H}\omega)\\
				&=L^*d(H\lrcorner\omega)+L^*(H\lrcorner d\omega)\\
				&=dh+L^*(H\lrcorner d\omega)
		\end{aligned}
	\end{equation}
	
	Both terms are zero, since $dh=0$ and $d\omega=0$.
\end{proof}

\begin{lemma}[\cite{dazord81}, \cite{morvan81}]
	If $L:\Sigma\rightarrow(X,\omega,J)$ is a Lagrangian submanifold of a Calabi-Yau manifold, then there is a smooth function $\beta:\Sigma\rightarrow S^1$, called the {\em Lagrangian angle}, with $h=d\beta$.
\end{lemma}

\begin{remark}
	The Lagrangian angle can be defined by the relation
		\begin{equation}
			L^*\Re(\Omega)=e^{i\beta}\dvol
		\end{equation}
	where $\Omega$ is the unit holomorphic $(m,0)$ form of the manifold $(X,\omega, J)$ and $\dvol$ is the volume element of $(\Sigma,g)$.
\end{remark}

\begin{definition}
	The class $[\lambda]\in H^1(\Sigma)$ is called the {\em period} or {\em Liouville class} of the immersion $L$. $[h]\in H^1(\Sigma)$ is the {\em Maslov class} of the immersion $L$. 
	
	If $[h]=0$, or equivalently if $\beta$ is a real-valued function, we say $L$ has {\em zero Maslov class}.
	
	If $[\lambda]=0$, or equivalently if $\lambda=d\phi$ for some smooth real-valued $\phi$, we say $L$ is {\em exact}.
\end{definition}

\begin{remark}
	If $X=\mathbf{C}^m$,  for any $x_0=(p_0^1,q_0^1,\ldots,p_0^m,q_0^m)\in\mathbf{C}^m$, we can take a Liouville form on $\mathbf{C}^m$ to be 
		\begin{equation}
			\eta_{x_0}=\sum_{i=1}^m (p^i-p_0^i)dq^i-(q^i-q_0^i)dp^i
		\end{equation}
	This choice of $\eta$ above shows that $\lambda$ is the one-form dual to the vector $L^*(-J(x-x_0)^\perp)$, where $(x-x_0)^\perp$ is the projection of the position vector $x-x_0$ to the normal bundle of $L$.
\end{remark}

\begin{remark}
Suppose $H^1(X)=0$. If $\eta_1$ and $\eta_2$ are two different Liouville forms on $(X,\omega,J)$ inducing $\lambda_1$ and $\lambda_2$ on $\Sigma$, then $d(\eta_1-\eta_2)=\omega-\omega=0$, so $\eta_1-\eta_2=df$ for some function $f$. 
	\begin{equation}
		\int_\gamma\lambda_1-\lambda_2 = \int_{L\circ\gamma}\eta_1-\eta_2=\int_{L\circ\gamma}df=0
	\end{equation}
Therefore $[\lambda]\in H^1(\Sigma)$ is independent of the choice of $\eta$. In particular, this is true when $(X,\omega,J)$ is the standard $\mathbf{C}^m$.
\end{remark}

We will use the Maslov class and period to study the singularities of the flow. We begin by recalling the following computations of Smoczyk \cite{smoczyk02}:

\begin{lemma}[\cite{smoczyk02}]
	The Maslov form and Liouville form evolve according to
	\begin{equation}\label{evol.forms}
		\begin{aligned}
			\frac{\partial}{\partial t}h&=dd^*h\\
			\frac{\partial}{\partial t}\lambda&=dd^*\lambda - 2h
		\end{aligned}
	\end{equation}
	where $d^*$ is the negative adjoint to $d$. In particular,
	\begin{equation}\label{evol.classes}
		\begin{aligned}
			\frac{\partial}{\partial t}[h]&=0\\
			\frac{\partial}{\partial t}[\lambda]&=-2[h]
		\end{aligned}
	\end{equation}
\end{lemma}

We will also use the scaling properties of $\lambda$ and $h$:
\begin{lemma}\label{scaling.lambda.h}
	Let $L:\Sigma\rightarrow (X,\omega,J)$ be a Lagrangian submanifold and $\tilde{L}:\Sigma\rightarrow(X,\alpha^2\omega,J)$ be its $\alpha$-rescale. Then $\tilde{\lambda}=\alpha^2\lambda$ and $\tilde{h}=h$, where $\tilde{\lambda}$ and $\tilde{h}$ are the Liouville and Maslov forms of $\tilde{L}$.
\end{lemma}
\begin{proof}
	If $\eta$ is a Liouville form for $(X,\omega,J)$, then $\alpha^2\eta$ is a Liouville form for $(X,\alpha^2\omega,J)$. So 
		\begin{equation}
			\tilde{\lambda}=L^*(\alpha^2\eta)=\alpha^2L^*\eta=\alpha^2\lambda
		\end{equation}
	
	To see how $h$ scales, note $\tilde{h}_{ijk}=-\alpha^2\omega(D_{\partial_i}\partial_j,\partial_k)=\alpha^2h_{ijk}$. Then
		\begin{equation}
			\begin{aligned}
				\tilde{h}=\tilde{H}_idx^i&=\tilde{g}^{jk}\tilde{h}_{ijk}dx^i\\
					&=\alpha^{-2}g^{jk}\alpha^2h_{ijk}dx^i=g^{jk}h_{ijk}dx^i=H_idx^i=h
			\end{aligned}
		\end{equation}
\end{proof}

Note that from Lemma \ref{scaling.lambda.h} it follows that the Lagrangian angle $\beta$ is scale-invariant and the primitive $\phi$ of the Liouville form has scaling degree 2, since the $d$ operator is scale-invariant.

\subsection{Rescaling Procedures}
The results and constructions in this section apply to mean curvature flow in general codimension in Euclidean space. For our purposes, we are concerned with  Lagrangian mean curvature flow in $\mathbf{C}^m$.

\subsubsection{The Tangent Flow}
	\begin{definition}
		If $\Sigma(t)$ has a singularity at time $T$ and $x_0\in \mathbf{C}^m$, and given $\displaystyle \alpha_j\rightarrow\infty$. The tangent flow of $\Sigma(t)$ is the limit $\tilde{\Sigma}_\infty(t)$ of the sequence of rescalings
		\begin{equation}
			\tilde{\Sigma}_j(t)=\alpha_j\left(\Sigma(T+\alpha_j^{-2}t)-x_0\right)
		\end{equation}
		where the limit is in the sense of Brakke flows, i.e.~one-parameter families of varifolds moving by their mean curvature.
	\end{definition}
	
	Here the scaling sequence $\alpha_j$ is arbitrary; of most interest to us will be the choices, if $t_j\rightarrow T$,  $\alpha_j=\left(T-t_j\right)^{-\frac{1}{2}}$ and $\displaystyle \alpha_j=\sup_{t\leq t_j}\lvert \II\rvert$.
	
	\begin{theorem}[Huisken]\label{TF.SS}
		Every tangent flow is a self-shrinker, i.e.
			\begin{equation}
				\tilde{\Sigma}^\infty(t)=\sqrt{-t}\tilde{\Sigma}^\infty(-1)
			\end{equation}
	\end{theorem}

\subsubsection{The Smooth Blow Up}

	\begin{theorem}\label{SBU}
		For any compact mean curvature flow $L:\Sigma\times[0,T)\rightarrow\mathbf{C}^m$, there is a sequence of rescalings $L_j(t)$ which converge, in $C^k$ for any $k$, on compact subsets of space and time to some smooth mean curvature flow
			\begin{equation*}
				L_\infty(t):\Sigma_\infty\times(-\infty,C)\rightarrow \mathbf{C}^m
			\end{equation*}
		where $C=\sup \lvert\II\rvert^2(T-t)$.
		
		This $L_\infty$ has achieves its space-time maximum of $\lvert\II\rvert^2$ at the origin in $\mathbf{C}^m$ and time $0$.
	\end{theorem}
	\begin{remark}
		The proof of Theorem \ref{SBU}, including the point-picking in the type II case, is quite standard. We sketch it here for completeness. A very clear exposition of the choices involved can be found in Chapter 8 of \cite{bensbook}.
		
	\end{remark}
	\begin{proof}
		Since we know that $\sup \lvert\II(p,t)\rvert=\infty$, we can select a sequence of points $(p_j,t_j)$ with $t_j\rightarrow T$ and $Q_j:=\lvert\II(p_j,t_j)\rvert=\max_{t\leq t_j}\lvert\II\rvert$. Moreover, in the type II case, we can select such a sequence with the property that $Q_j^2(T-t_j)=\sup_{t\leq t_j}\lvert II(p,t)\rvert^2(T-t)$.

	We define, for a such sequence of space-time points $(p_j,t_j)$ along which $\lvert \II\rvert$ is blowing up, rescalings
			\begin{equation}
				L_j(t)=Q_j\left(L(t_j+Q_j^{-2}t)-x_j\right)
			\end{equation}
		where $x_j=L(t_j,p_j)$.
		
		This construction guarantees that each rescaled Lagrangian $L_j$ has its second fundamental form $\lvert\II_j(\cdot,t)\rvert$ bounded by $1$ for $t\leq 0$; in fact that $\lvert\II_j\rvert$ is bounded on compact subsets of $[-Q_j^2T,Q_j^2(T-t_j))$. The mean curvature flow equation allows us to conclude that the covariant derivatives of each $\II_j$ are bounded on the same set. Thus the $\Sigma_j$ converge geometrically to some $L_\infty:\Sigma_\infty\times(-\infty,C)\rightarrow\mathbf{C}^m$, where $C=\limsup Q_j^2(T-t_j)$.
	\end{proof}
	
	\begin{definition}
		We will refer to $L_\infty(t)$ as the smooth blow-up of the singularity.
	\end{definition}
	\begin{remark}
		The domain manifold $ \Sigma_\infty$ of the smooth blow-up will in general be different from the original domain manifold $\Sigma$; in particular it may be noncompact.
	\end{remark}

	\begin{theorem}\label{SBU.SS}
		If the singularity of $L:\Sigma\times[0,T)\rightarrow\mathbf{C}^m$ is of type I, then the smooth blow-up $L_\infty:\Sigma_\infty\times(-\infty,C)\rightarrow\mathbf{C}^m$ is a smooth self-shrinker with singularity at $(x_\infty, C)$, where ${\displaystyle C=\limsup_j Q_j^2(T-t_j)}$.
	\end{theorem}

	\begin{remark}
		The proof of Theorem \ref{SBU.SS} is largely the same as the proof of Theorem \ref{TF.SS}, with the type I hypothesis supplying the necessary control on the Huisken functional $\Theta$. 
		
		Since the proofs are so similar and the proof of Theorem \ref{TF.SS} is easily found elsewhere, we omit them. The proof of Theorem \ref{BDSS} below is very much in the same vein.
	\end{remark}

\subsubsection{Blowing down the smooth blow-up}
	To obtain the smooth blow-up, we must rescale both at a different space-time center and by a  larger factor than in the tangent flow. In this section we will show that, in fact, the smooth blow-up of a type II singularity can be ``blown down" to recover a singularity model very similar to the tangent flow construction.
	
	\begin{definition}
		If $\Sigma(t)$, $t\in[0,T)$ is a mean curvature flow with singularity at time $T$ and $\Sigma_\infty(t)$ is a smooth blowup, and $\epsilon_j\rightarrow 0$, the {\em blow-down} of $\Sigma_\infty$ is the limit $\Sigma_\infty^\infty$ of the rescalings
			\begin{equation*}
				\Sigma_\infty^j=\epsilon_j\Sigma_\infty\left(\epsilon_j^{-2}t\right)
			\end{equation*}
		If we pick $\epsilon_j$ by $\epsilon_j^{-2}=Q_j(T-t_j)^{\frac{1}{2}}$ then we call the corresponding $\Sigma_\infty^\infty$ a {\em recovery blow-down}.
	\end{definition}

	\begin{remark}
		Since the original singularity is of type II, $\epsilon_j\rightarrow 0$ and $\epsilon_j^{-2}\rightarrow \infty$, so the blow-down captures the behavior of $\Sigma_\infty$ at infinite backward time; it is a backward limit. We will argue below that in fact $\Sigma_\infty$ captures the singular behavior of the original flow $\Sigma(t)$ in some sense.
		
		In the type I case, when $\Sigma_\infty(t)$ develops a singularity at $x\in\mathbf{C}^m$ at time $C<\infty$, we rescale about $(x,C)$ by any sequence $\epsilon_j\rightarrow 0$.
	\end{remark}
	
	\begin{theorem}\label{BDSS}
		The blow-down $\Sigma_\infty^\infty$ of the smooth blow-up of a compact mean curvature flow is a self-shrinker.
	\end{theorem}
		The proof of Theorem \ref{BDSS} is essentially standard, using the Huisken functional
			\begin{equation}
				\begin{aligned}
				\Theta_{(x_0,T)}(t)&=\int_{\Sigma(t)}\left(4\pi(T-t)\right)^{-\frac{m}{2}}\exp\left(-\frac{\lvert x-x_0\rvert^2}{4(T-t)}\right)\dHaus^m\\
				&=\int_{\Sigma(t)}\theta_{(x_0,T)}(x,t)\dHaus^m
				\end{aligned}
			\end{equation}
		which evolves along the mean curvature flow by
			\begin{equation}
				\frac{d}{dt}\Theta_{(x_0,T)}(t)=-\int_{\Sigma(t)}\left\lvert H+\frac{(x-x_0)^\perp}{2(T-t)}\right\rvert^2\theta_{(x_0,T)}(x,t)\dHaus^m
			\end{equation}
	and hence is monotone along the flow. We define $\Theta(x_0,T)$ to be ${\displaystyle \lim_{t\rightarrow T} \Theta_{(x_0,T)}(t)}$.

	We will also need the corresponding entropy functional, which is also monotone along the flow \cite{colding&minicozzi09}:
		\begin{equation}
			\Lambda(t)=\sup\Theta_{(x_0,T)}(t)
		\end{equation}
	where the supremum is taken over all $x_0\in\mathbf{C}^m$ and all $T>t$. We define ${\displaystyle \Lambda(T)=\lim_{t\rightarrow T} \Lambda(t)}$.

	It follows directly from the scale-invariance of $\Theta_{(x_0,T)}$ under parabolic rescaling about $(x_0,T)$ that
	\begin{proposition}\label{entropybound}
		The entropy $\Lambda_\infty(t)$ of $L_\infty(t)$ is bounded by $\Lambda(T)$.
	\end{proposition}
	
	In particular, the Huisken's functional 
		\begin{equation}
			\Theta^\infty(t)=\int_{\Sigma_\infty(t)}\theta_{(0,0)}(x,t)\dHaus^m
		\end{equation}
	is finite and uniformly bounded for all $t<0$.
	
	\begin{proof}[Proof of Theorem \ref{BDSS}]
		For compact $K\subset\mathbf{C}^m$ and $a<b<0$, we consider
		\begin{equation}\label{BDSS.ineq}
			\begin{aligned}
				\int_a^b\int_{\Sigma_\infty^j(t)\cap K}\left\lvert H-\frac{x^\perp}{2t}\right\rvert^2\theta_{(0,0)}(x,t)\dHaus^mdt\\
				=\int_{\epsilon_j^{-2}a}^{\epsilon_j^{-2}b}\int_{\Sigma_\infty(t)\cap \epsilon_jK}\left\lvert H-\frac{x^\perp}{2t}\right\rvert^2\theta_{(0,0)}(x,t)\dHaus^mdt\\
				\leq \Theta^\infty(\epsilon_j^{-2}a)-\Theta^\infty(\epsilon_j^{-2}b)
			\end{aligned}
		\end{equation}
		which goes to $0$ as $j\rightarrow\infty$. Thus the quantity $\left\lvert H-\frac{x^\perp}{2t}\right\rvert^2$ vanishes in the limit, which forces $L_\infty^\infty(t)$ to be a self-shrinking flow, i.e.
		\begin{equation}
			L_\infty^\infty(t)=\sqrt{-t}L_\infty^\infty(-1)
		\end{equation}
	\end{proof}

\section{Type-II SBUs are Exact and Zero-Maslov}
In this section we investigate the behavior of the classes $[\lambda]$ and $[h]$ under the smooth blow-up procedure. To this end, suppose $\gamma_\infty$ is a 1-cycle in $\Sigma_\infty$. Since $\gamma_\infty$ is compact, by geometric convergence there is a corresponding $\gamma_j$ in $\Sigma_j=\Sigma$ for large enough $j$. We have:

	\begin{align}
		\lambda_\infty(t).\gamma_\infty=&\lim_j \lambda_j(t).\gamma_j\\
		h_\infty(t).\gamma_\infty=&\lim_j h_j(t).\gamma_j
	\end{align}
where
	\begin{equation}
		\lambda_j(t)=\left(L_j(t)\right)^*\lambda_{x_j}
	\end{equation}

Moreover, by geometric convergence the $\gamma_j$ are homologous for large enough $j$. So there is some 1-cycle $\gamma\subset \Sigma$ with

	\begin{align}
		\lambda_\infty(t).\gamma_\infty&=\lim_j \lambda_j(t).\gamma\\
		h_\infty(t).\gamma_\infty&=\lim_j h_j(t).\gamma
	\end{align}

A simple computation shows that 
	\begin{align}
		\lambda_j(t).\gamma&=Q_j^2\lambda\left(t_j+Q_j^{-2}t\right).\gamma\\
		h_j(t).\gamma&=h\left(t_j+Q_j^{-2}t\right).\gamma
	\end{align}
	
We will combine these rescalings with the evolution equations (\ref{evol.forms}), to obtain the following theorem:

	\begin{namedtheorem}[Theorem A]\label{cohomological}
		Let $L_\infty:\Sigma_\infty\times(-\infty,\infty)\rightarrow\mathbf{C}^m$ be the smooth blow-up of a compact Lagrangian mean curvature flow $L:\Sigma\times[0,T)\rightarrow\mathbf{C}^m$. Writing $h$, $\lambda$ for the Maslov and Liouville forms of $L$ and $h_\infty$, $\lambda_\infty$ for the Maslov and Liouville forms of $L_\infty$, we have:
		\begin{enumerate}
			\item If the singularity is type I, then  $[\lambda_\infty(t)]=2(C-t)[h_\infty]$, where $C=\limsup_j Q_j^2(T-t_j)$. Any cycle $\gamma\subset\Sigma$ for which $[\gamma]$ survives to a class $[\gamma_\infty]\in H_1(\Sigma_\infty;\mathbf{R})$ has $\lambda(0).\gamma=2Ch(0).\gamma$.
			\item If the singularity is type II, then $[\lambda_\infty]=0$ and $[h_\infty]=0$. Moreover, any cycle $\gamma\subset\Sigma$ for which $[\gamma]$ survives to a class $[\gamma_\infty]\in H_1(\Sigma_\infty;\mathbf{R})$ has $\lambda.\gamma=0$ and $h.\gamma=0$.
		\end{enumerate}
	\end{namedtheorem}
	
	\begin{remark}
		Actually the first clause follows from Theorem \ref{SBU.SS}, since in the type I case, the SBU is a self-shrinker, so we have $\lambda_\infty(t)=2(C-t)h_\infty(t)$ as forms, not merely at the level of cohomology.
	\end{remark}

	\begin{proof}
		If $\gamma_\infty$ is a 1-cycle in $\Sigma_\infty$ coming from $\gamma\subset\Sigma$, we have
			\begin{equation}\label{rescale.comp}
				\begin{aligned}
					\lambda_\infty(t).\gamma_\infty&=\lim_jQ_j^2\lambda\left(t_j+Q_j^{-2}t\right).\gamma\\
					&=\lim_j\left(Q_j^2\left[\lambda(0).\gamma-2\left(t_j+Q_j^{-2}t\right)h(0).\gamma\right]\right)\\
					&=\lim_j\left(Q_j^2\left[\lambda(0).\gamma-2t_jh(0).\gamma\right]\right)-2th(0).\gamma
				\end{aligned}
			\end{equation}
		To make the left-hand side of (\ref{rescale.comp}) finite, it must be the case that
			\begin{equation}
				\lim_j\left[\lambda(0).\gamma-2t_jh(0).\gamma\right]=0
			\end{equation}
		i.e., that
			\begin{equation}
				\lambda(0).\gamma=2Th(0).\gamma
			\end{equation}
		so that $\gamma\in\ker\left[\lambda(0)-2Th(0)\right]$. We continue (\ref{rescale.comp}):
			\begin{equation}\label{rescale.comp.cont}
				\begin{aligned}
					\lambda_\infty(t).\gamma_\infty&=\lim_j2Q_j^2\left(T-t_j\right)h(0).\gamma-2th(0).\gamma
				\end{aligned}
			\end{equation}
		
		Again, for this limit to be finite, we must have either that $\lim_jQ_j^2\left(T-t_j\right)$ is finite, i.e. the singularity is type I, or that $h(0).\gamma=0$.
		
		In the former case, the computations (\ref{rescale.comp}) and (\ref{rescale.comp.cont}) give the first clause of Theorem A.
		
		Similarly in the latter case, $h(0).\gamma=0$  gives $\lambda_\infty(t).\gamma_\infty=0$ by (\ref{rescale.comp.cont}). Since the Maslov class is scale-invariant, $h(0).\gamma=0$ implies $h_\infty(t).\gamma_\infty=0$.
	\end{proof}

\section{A Uniform Bound for the Lagrangian Angle}
	Theorem A shows that, for any singularity which occurs at a cohomologically-incorrect time, the Lagrangian angle $\beta$ lifts to a real-valued function near the singularity. We would like to show that this real-valued lift is uniformly bounded.
	
	In this section, we prove a kind of maximum principle for circle-valued functions on a closed manifold which satisfy the time-dependent heat equation.
	
	First notice that while the notions of ``maximum" and ``minimum" do not make sense for a circle-valued function, the notion of oscillation does. Given a circle-valued function $\beta:\Sigma\rightarrow S^1$ on a closed manifold, by elementary homotopy theory there is a lift $\tilde{\beta}:\tilde{\Sigma}\rightarrow \mathbf{R}$, where $\tilde{\Sigma}$ is the universal cover of $\Sigma$. 
	\begin{definition}
		Let $D$ be a fundamental domain in $\tilde{\Sigma}$ and $\tilde{\beta}$ be a lift of $\beta$. The oscillation of $\beta$ is
			\begin{equation}
				\osc \beta = \osc_{D}\tilde{\beta}=\max_{D}\tilde{\beta} - \min_D \tilde{\beta}
			\end{equation}
	\end{definition}
	
	\begin{proposition}
		$\osc \beta$ is well-defined independent of choice of lift.
	\end{proposition}
	\begin{proof}
		If $\tilde{\beta}_1$ and $\tilde{\beta}_2$ are two lifts of $f$, then we have for each $x\in \tilde{\Sigma}$ $e^{i\tilde{\beta}_1(x)}=e^{i\tilde{\beta}_2(x)}=\beta(\pi(x))$. So $e^{i(\tilde{\beta}_2-\tilde{\beta}_1)}=1$, so $\tilde{\beta}_2-\tilde{\beta}_1$ is a continuous lift of the constant function $1$. The only such lifts are multiples of $2\pi$. Therefore $\tilde{\beta}_2=\tilde{\beta}_1+2k\pi$ for some $k$. Then
			\begin{equation}
				\max_D\tilde{\beta}_2-\min_D\tilde{\beta}_2=\max_D\tilde{\beta}_1 + 2k\pi - \min_D\tilde{\beta}_1-2k\pi=\max_D\tilde{\beta}_1-\min_D\tilde{\beta}_1
			\end{equation}
	\end{proof}
	
	\begin{remark}
		We also note that $\osc \beta$ is independent of the choice of fundamental domian $D$.
	\end{remark}

	\begin{theorem}[$S^1$-Maximum Principle]\label{maxprin.metric}
		Suppose $\beta:\Sigma\rightarrow S^1$ satisfies the time-dependent heat equation $(\frac{\partial}{\partial t}-\Delta)\beta=0$. Then $\osc \beta(t)$ is bounded independent of time.
	\end{theorem}
	The idea of the proof is to compare the function $\beta$ to a harmonic function in the same topological class.	
	\begin{proof}	
			Let $\eta(t)$ be the harmonic representative of the class $[d\beta]$ with respect to the metric at time $t$. Then $\eta-d\beta=df$ for some function $f$. Now
				\begin{equation}\label{comp.heat.metric}
					\begin{aligned}
					\left(\frac{\partial}{\partial t}-\Delta\right)df&=\left(\frac{\partial}{\partial t}-\Delta\right)(\eta-d\beta)\\
					&=\frac{\partial}{\partial t}\eta-(d^*d+dd^*)\eta-d\frac{\partial}{\partial t}\beta + (d^*d+dd^*)d\beta\\
					&=\frac{\partial}{\partial t}\eta-d\left((\frac{d}{dt}-\Delta)\beta\right)=\frac{\partial}{\partial t}\eta
					\end{aligned}
				\end{equation}
			since $\beta$ solves the heat equation.
			
			As above, $\eta=dh$ for some $S^1$-valued function $h$, which we can take to be harmonic at each time. We have
			
			\begin{equation}
				\left(\frac{\partial}{\partial t}-\Delta\right)f=\frac{\partial}{\partial t}h
			\end{equation}
			
			So by the ordinary scalar maximum principle with reaction term, for each $x\in \Sigma$ we have
			\begin{equation}
				f(x,t)\leq \max f(0)+\int_0^t\frac{\partial }{\partial t}h(x,\tau)d\tau
			\end{equation}
			Lifting to the fundamental domain, we have 
			\begin{equation}
				f(x,t)\leq \max f(0)+\tilde{h}(x,t)-\tilde{h}(x,0)
			\end{equation}
			so that
			 \begin{equation}
				\max f(t)\leq \max f(0)+\max_D \tilde{h}(t)-\min_D\tilde{h}(0)
			\end{equation} 
			and similarly 
			 \begin{equation}
				\min f(t)\geq \min f(0)+\min_D \tilde{h}(t)-\max_D\tilde{h}(0)
			\end{equation} 
			so that
			\begin{equation}
				\osc f(t)\leq \osc f(0)+\osc h(t)-\osc h(0)\leq\osc f(0)+\osc h(t)
			\end{equation}
			
			So we have, using the definition of $f$,
			\begin{equation}
				\osc\beta(t)\leq \osc h(t)+\osc f(t)\leq \osc f(0)+2\osc h(t)
			\end{equation}
			
			Now we want to give a bound for $\osc h(t)$. 
		\begin{lemma}\label{osc.harmonic} 
		 	$\osc h(t)$ is bounded topologically, i.e.~in terms of $[h]\in H^1(\Sigma)$.
		\end{lemma}
		\begin{proof}
		If $D$ is a fundamental domain of the universal cover of $\Sigma$, with $\partial D$ having faces $F_1,\ldots F_{2p}$, we can choose $[\gamma_1],\ldots,[\gamma_p]$  generating $\pi_1(M)$ so that the deck transformation $T_k$ corresponding to $[\gamma_k]$ takes $F_{2k-1}$ to $F_{2k}$. We define $g_k:\partial D\rightarrow [0,1]$ as follows. $g_k\equiv 0$ on $F_{2k-1}$ and $g_k\equiv 1$ on $F_{2k}$.  On other faces, define $g_k$ so that $g_k|_{F_{2j}}\circ T_j=g_k|_{F_{2j-1}}$.
							
		Let $h_k$ be the harmonic function on $D$ which solves the Dirichlet problem $h_k|_{\partial D}=g_k$. Since each $h_k$ is harmonic, we have $\osc_D h_k=\osc_{\partial D} h_k=1$. Moreover, by construction $h_k$ actually defines a circle-valued function on $\Sigma$. So we can consider $[dh_k]\in H^1(\Sigma;\mathbf{R})$. Under the homomorphism
			\begin{equation}
				\pi_1(\Sigma)\rightarrow H_1(\Sigma;\mathbf{Z})\rightarrow H_1(\Sigma;\mathbf{R})\cong H^1(\Sigma;\mathbf{R})
			\end{equation} 
		(where the first map is the abelianization map, the second is tensor product with $\mathbf{R}$, and the last is the Poincar\'e duality isomorphism), we have that $[\gamma_k]\mapsto [dh_k]$. So the $[dh_k]$ span $H^1(\Sigma;\mathbf{R})$.
		
		We can therefore write $[d\beta]=[dh]=\sum \alpha_k[dh_k]$, and since the forms $dh,dh_k$ are harmonic, the Hodge theorem gives $dh=\sum \alpha_k\ dh_k$. Up to addition of a constant, therefore, we have $\tilde{h}=\sum_k\alpha_kh_k$. On the other hand, each $h_k$ has oscillation 1, so $\osc h\leq \sum_k\lvert \alpha_k\rvert$.
		
		 The coefficients $\alpha_k$ may not be unique, since $[dh_k]$ span but may not be independent; but we can take the infimum over all choices of $\alpha_k$, which is determined by the classes $[dh_k], [d\beta]\in H^1(\Sigma;\mathbf{R})$, which depend in turn only on the choices of fundamental domain $D$ and generating loops $[\gamma_1],\ldots,[\gamma_p]$. 
		 \end{proof}
			
			Thus the oscillation of $h(t)$ is bounded topologically, hence independent of $t$. The bound on $\osc\beta$ follows. \end{proof}

	\begin{theorem}\label{betabound}
		Suppose the smooth blow-up $\Sigma_\infty$ is zero-Maslov. Then its Lagrangian angle is uniformly bounded on any time interval of the form $(-\infty, t_0]$.
	\end{theorem}
	The proof of this result involves the use of the fact that the Lagrangian angle satisfies the heat equation. We use a conjugate heat kernel $\Psi_{p_0,t_0}(p,t)$ on $\Sigma_\infty$. Since the smooth blow-up is in general noncompact, the existence of a backward head kernel must be justified. Perelman stated \cite{perelman02} and S.~Zhang gave a complete proof \cite{szhang10} that under geometric convergence, if $\Psi_k$ is the conjugate heat kernel for a sequence of geometrically converging Ricci flows $(M_k,g_k(t),p_k)\rightarrow (M_\infty,g_\infty(t),p_\infty)$, that is if:
	\begin{gather}
		\left\{
		\begin{aligned}
			\left(-\frac{\partial}{\partial t}-\Delta_{g_k(t)}+R(g_k(t))\right)u_k=&0\\
			\lim_{t\nearrow t_0} u_k(t)=&\delta_{p_k}
		\end{aligned}
		\right.
	\end{gather}
	then the $\Psi_k$ subconverge on compact domains to a conjugate heat kernel $\Psi_\infty$ on $(M_\infty, g_\infty(t),p_\infty)$ with center $(p_\infty,t_0)$. Zhang's proof is stated for the Ricci flow, but it relies only upon estimates of Chau-Tam-Yu which hold for evolving metrics in general \cite{chau&tam&yu11}. The convergence given by Theorem \ref{SBU} in particular implies Cheeger-Gromov convergence for the one-parameter families $(\Sigma_k,g_k(t),p_k)$. We state the theorem in detail:
	
\begin{theorem}[Zhang's theorem for MCF]\label{convergence.kernels}
	Suppose $F_k:(M_k,p_k)\times[\alpha,\omega]\rightarrow (N_k,h_k,x_k)$ are a sequence of mean curvature flows converging geometrically to $F_\infty:(M_\infty,p_\infty)\times[\alpha,\omega]\rightarrow (N_\infty,h_\infty)$. Let $t_0\in[\alpha,\omega]$.  Let $g_k(t)=F_k^*(t)h_k$ and $\Psi_k$ be conjugate heat kernel at $(p_k,t_0)$, i.e. solve minimally:
		\begin{gather}\label{conjugate.heat.MCF}
			\left\{
			\begin{aligned}
				\left(-\frac{\partial}{\partial t}-\Delta_{g_k(t)} + \lvert H_k\rvert^2\right)\Psi=&0\\
				\lim_{t\nearrow t_0} \Psi(t)=&\delta_{p_k}
			\end{aligned}
			\right.
		\end{gather}
	If $\phi_k:(W_k,p_\infty)\rightarrow (M_k,p_k)$ are the convergence embeddings, then $\tilde{\Psi}_k=\Psi_k\circ\phi_k$ converge on compact domains to $\Psi_\infty$, to the conjugate heat kernel at $(p_\infty,t_0)$.
\end{theorem}

Chau-Tam-Yu's estimates guarantee that, at each time, the conjugate heat kernel decays fast enough at (spatial) infinity to justify integration by parts against any function which at most polynomial (spatial) growth. In particular:

\begin{proposition}\label{duh}
	Let $f\in C^\infty(M\times[\alpha,\omega])$ be a smooth function on a mean curvature flow $F_t:M\times[\alpha,\omega]\rightarrow \mathbf{R}^{m+n}$ with at most polynomial growth at each time. If $u$ is a conjugate heat kernel, then the $u$-weighted integral of $f$ evolves according to
		\begin{equation*}
			\frac{d}{dt}\int_{M}f(t)u(t)\dvol(t)=\int_{M}(\frac{\partial}{\partial t}-\Delta)f(t)u(t)\dvol(t)
		\end{equation*}
\end{proposition}
	
	\begin{proof}[Proof of Theorem \ref{betabound}]
		Let $\beta$ be the Lagrangian angle of the original flow $\Sigma(t)$. By the maximum principle above, $\osc\beta$ is uniformly bounded in time. Since the Lagrangian angle is scale-invariant, this means that $\osc\beta_\infty$ is also uniformly bounded in time along the flow $\Sigma_\infty$.

		Pick any $p_0\in \Sigma_\infty$. Let $\Psi_{(p,t_0)}(t,y)$ be the conjugate heat kernel along the flow, centered at $p_0,t_0$. Since $\beta_\infty$ satisfies the heat equation, we have by Proposition \ref{duh} that  
			\begin{equation}
				\int_{\Sigma_\infty(t)}\beta_\infty(t) \Psi_{(p_0,t_0)}(t)\dHaus^m\equiv \beta_\infty(p_0,t_0)
			\end{equation}
		Therefore we define the function (for $t<t_0$), $\overline\beta_\infty(t)=\beta_\infty(t)-\beta_\infty(p_0,t_0)$, which has $\int\overline\beta_\infty(t)\Psi_{(p_0,t_0)}(t,y)\dvol(t,y)\equiv 0$. In particular, for each $t$, $\min\overline{\beta}_\infty(t)<0<\max\overline{\beta}_\infty(t)$.
		
		But then $\max \overline\beta_\infty(t)<\osc\overline\beta_\infty(t)=\osc\beta_\infty(t)$, which is uniformly bounded. Similarly, $\min\overline\beta_\infty$ is uniformly bounded. So $\overline\beta_\infty$ is bounded uniformly in time. But $\overline{\beta}_\infty$ differs from $\beta_\infty$ by a constant, so $\beta_\infty$ is uniformly bounded on $(-\infty,t_0]$.
	\end{proof}

\section{Blowing Down Zero-Maslov Solutions}
	In this section we will establish the following theorem
		\begin{theorem}\label{zeroMaslovSBU}
			Suppose the smooth blow-up is zero-Maslov. Then any blow-down is a union of special Lagrangian cones $L_1,\ldots,L_N$. Moreover, the convergence respects the Lagrangian angle in the sense that if $\theta_1,\ldots,\theta_N$ are the angles of the special Lagrangian cones, we have
			\begin{equation*}
				\lim_k\int_{\Sigma_\infty^k}f(\exp(i\beta_\infty))\phi \dHaus^m=\sum_j f(\exp(i\theta_j))m_j\mu_j(\phi)
			\end{equation*}
			for every continuous function $f: S^1\rightarrow \mathbf{R}$ and every compactly supported continuous function $\phi$ on $\mathbf{C}^m$, where  $m_j$ and $\mu_j$ are the multiplicity and measure associated to the special Lagrangian cone $L_j$.
		\end{theorem}
	The proof of this theorem proceeds along the same lines as Neves' Theorem \cite{neves&tian07}
	
	\begin{proof}
	By Theorem \ref{betabound}, the Lagrangian angle $\beta_\infty$ is uniformly bounded on $\Sigma_\infty\times (-\infty,0]$. Therefore, using Proposition \ref{entropybound} the quantity
		\begin{equation}
			\mathcal{B}(t)=\int_{\Sigma_\infty(t)}\beta_\infty^2\theta_{(0,0)}(t)\dHaus^m
		\end{equation}
	is uniformly bounded for $t\in(-\infty,0]$.
	
	Because the smooth blow-up $L_\infty$ is noncompact, we use the Huisken's kernel $\theta_{(0,0)}$ to weight all integral estimates. It satisfies the following
	\begin{proposition}[kernel property of $\theta$]\label{CHK}
		If $\Sigma(t)$ is a mean curvature flow and $f$ is a smooth scalar function, then
			\begin{equation*}
				\begin{aligned}
					\frac{d}{dt}\int_{\Sigma(t)}f\theta_{(x_0,T)}\dHaus^m=&-\int_{\Sigma(t)}\left\lvert H +\frac{1}{2(T-t)}(x-x_0)^\perp\right\rvert^2f\theta_{(x_0,T)}\dHaus^m\\
					&+\int_{\Sigma(t)}\left(\frac{\partial}{\partial t}f-\Delta f\right)\theta_{(x_0,T)}\dHaus^m
				\end{aligned}
			\end{equation*}
	\end{proposition}
	
	Since $\beta_\infty$ satisfies the heat equation,
		\begin{equation}
			\frac{\partial}{\partial t}\beta_\infty^2-\Delta\beta_\infty^2=-2\lvert\nabla \beta_\infty\rvert^2=-2\lvert H\rvert^2
		\end{equation}
	and applying Proposition \ref{CHK}, we have
		\begin{equation}
			\begin{aligned}
			\frac{d}{dt}\mathcal{B}(t)=-2&\int_{\Sigma_\infty(t)}\lvert H\rvert^2\theta_{(0,0)}\dHaus^m\\
			-&\int_{\Sigma_\infty(t)}\beta_\infty^2\left\vert H-\frac{1}{2t}x^\perp\right\vert^2\theta_{(0,0)}\dHaus^m
			\end{aligned}
		\end{equation}
	Thus $\mathcal{B}(t)$ is monotone, and constant precisely along minimal self-shrinkers, i.e. minimal cones.
	
	On the other hand, for any $a<b<0$ and compact $K\subset\mathbf{C}^m$, we have (as in (\ref{BDSS.ineq}))
		\begin{equation}
			\begin{aligned}
				2\int_a^b\int_{\Sigma_\infty^j(t)\cap K}\lvert H\rvert^2\theta_{(0,0)}\dHaus^mdt&\leq \mathcal{B}(\epsilon_j^{-2}a)-\mathcal{B}(\epsilon_j^{-2}b)
			\end{aligned}
		\end{equation}
		
	Now taking the limit in $j$, since $\mathcal{B}(t)$ is monotone and bounded, both $\mathcal{B}(\epsilon_j^{-2}a)$ and $\mathcal{B}(\epsilon_j^{-2}b)$ approach the same limit. So 
		\begin{equation}
			\lim_j \int_a^b\int_{\Sigma_\infty^j(t)\cap K}\lvert H\rvert^2\theta_{(0,0)}\dHaus^mdt=0
		\end{equation}
	In particular, for almost all $t\in[a,b]$ we have
		\begin{equation}
			\lim_j\int_{\Sigma_\infty^j(t)\cap K}\lvert H\rvert^2\theta_{(0,0)}\dHaus^m=0
		\end{equation}. 
	By a similar argument to that in the proof of Proposition 5.1 of \cite{neveszeromaslov}, this shows that $\Sigma_\infty^\infty$ is a union of special Lagrangian cones, with convergence of the Lagrangian angle functions.\end{proof}

	Since by Theorem A every type II smooth blow-up is zero-Maslov, Theorem \ref{zeroMaslovSBU} immediately implies one of our main theorems:
	\begin{namedtheorem}[Theorem B]
		All type II singularities are modeled by special Lagrangian cones.
	\end{namedtheorem}

\section{Type I singular times}
	In this section we will establish a converse to the first clause of Theorem A, namely
	\begin{namedtheorem}[Theorem C]\label{my.toosoon}
		Suppose $L(0)$ is a compact Lagrangian with $\ker[h]=0$ and finitely-generated first homology, generated by $\gamma_1,\ldots,\gamma_k$. Then a type I singularity of the Lagrangian mean curvature flow starting from $L(0)$ can only occur at one of
		\begin{equation*}
			T_1=\frac{\int_{\gamma_1}\lambda}{2\int_{\gamma_1}h},\ldots,T_k=\frac{\int_{\gamma_k}\lambda}{2\int_{\gamma_k}h}
		\end{equation*}
	\end{namedtheorem}
	
	This theorem extends work of Groh-Schwarz-Smoczyk-Zehmisch and Neves from the monotone surface case to the general case.
	\begin{theorem}[\cite{groh&schwarz&smoczyk&zehmisch} \cite{nevesmonotone}]\label{smoczyk.toosoon}
		A singularity of an embedded normalized monotone Lagrangian mean curvature flow which occurs at time $T<\frac{1}{2}$ must be of type II.
	\end{theorem}

	\begin{proof}[Proof of Theorem C]
		Suppose $T$ is not one of the cohomologically prescribed times. Then by Theorem A, no cycle survives. Therefore the smooth blow-up $\Sigma_\infty$ is zero-Maslov. Thus as before, its blow-down must be a minimal cone.
		
		On the other hand, since the smooth blow-up of a type I singularity is self-similar, the blow-down is the same as the smooth blow-up. Therefore the smooth blow-up is a nonflat smooth minimal cone, which is a contradiction.
	\end{proof}

	Grayson \cite{grayson87} completely characterized the singularities of embedded curve-shortening flow in the plane: starting from an initial curve $\Gamma\subset \mathbf{R}^2$, which bounds an area $A$, the flow terminates at time $\frac{A}{2\pi}$ in a ``round point". In subsequent papers he considered the case of index-zero curves, showing that their singularities are all modeled by the ``grim reaper" curve $y=\log\cos x$. Consider the product $\Gamma_1\times\cdots\times \Gamma_m\subset\mathbf{C}^m$, where each $\Gamma_i$ is an embedded plane curve. Then Grayson's result says that a singularity will occur exactly at time $T=\min\frac{A_i}{2\pi}$, where $A_i$ is the area bounded by $\Gamma_i$; moreover the singularity will be modeled by a product whose factors are $\mathbf{R}$ or $S^1$, with the $S^1$ factors coming from those curves bounding the least area. Theorem C can be viewed as an extension of this result.

	In the case of a product of embedded curves $\Gamma_1\times\cdots\times\Gamma_m$, singularities happen because one factor collapses before the others. One might suspect this is the case in general:

	\begin{conjecture}
		In Theorem C, a type I singularity can only occur at the first of the candidates $T_i$.
	\end{conjecture}

\section{Minimal SBUs}
Recall the almost-calibrated clause of Neves' theorem on zero-Maslov singularities:
	\begin{reptheorem}{TF.zeromaslov}[\cite{neveszeromaslov}]
		If $\Sigma(0)$ is almost-calibrated, that is, if $\cos\beta\geq \epsilon>0$, then any tangent flow is a minimal cone with a single Lagrangian angle.
	\end{reptheorem}
	
	Neves posed the question of whether, in fact, the smooth blow-up of such a singularity is minimal. We have an affirmative answer:
	\begin{namedtheorem}[Theorem B', first clause]
		If the initial data are almost-calibrated, then any smooth blow-up is a smooth special Lagrangian.
	\end{namedtheorem}
	\begin{proof}
		Neves' argument gives a tangent flow with a single Lagrangian angle. Thus the Lagrangian angle $\beta^\infty$ must at each time asymptotically approach $\beta_0$ at spatial infinity. But $\beta_\infty$ also satisfies the heat equation, so it can have no spacetime interior local minima or local maxima.
		
		Thus we must have $\beta^\infty\equiv \beta_0$, so that $\lvert H_\infty\rvert=\lvert h_\infty\rvert=\lvert d\beta^\infty\rvert=0$.
	\end{proof}
	
	Similarly, in the dimension two case of a monotone flow with a premature singularity, Neves showed
	\begin{reptheorem}{TF.monotone}[\cite{nevesmonotone}]
		Suppose $m=2$, $\Sigma(0)$ is monotone $[\lambda]=C[h]$, and $T<\frac{C}{2}$. Then any tangent flow is a minimal cone with a single Lagrangian angle.
	\end{reptheorem}
	
	Theorem A says that in case of a premature singularity, the smooth blow-up is zero-Maslov. So we can use the above proof to obtain:
	\begin{namedtheorem}[Theorem B', second clause]
		Under the hypotheses of Theorem \ref{TF.monotone}, any smooth blow-up is a smooth special Lagrangian.
	\end{namedtheorem}

	\section{Some Questions}
	
	\subsection{Singularity type and embeddedness}
	In general, high-codimension mean curvature flow does not preserve embeddedness. In the Lagrangian case, we can relate embeddedness to $\lambda$ and $h$ as follows.

	Viterbo \cite{viterbo90} showed that, for any Lagrangian embedding $L:\Sigma^m\rightarrow\mathbf{C}^m$ of a closed manifold $\Sigma$ which admits a metric of negative curvature, there is at least one cycle $\gamma$ with (in our notation) $h.\gamma=2\pi$ and $\lambda.\gamma>0$. In particular this includes all closed surfaces of higher genus. Viterbo's proof also holds for $\Sigma=T^2$. 

		\begin{remark}
			Polterovich \cite{polterovich91} gave related results following Gromov's famous work. For a summary of how these results are related, see the introduction by Oh \cite{oh96}.
		\end{remark}

	\begin{definition}
		We call any cycle with $h.\gamma=2\pi$ and $\lambda.\gamma>0$ a {\em Viterbo cycle}. Given a Lagrangian surface evolving by mean curvature with embedded initial data, define $\Lambda^+=\max \lambda(0).\gamma$ and $\Lambda^-=\min\lambda(0).\gamma$, where the maximum and minimum are taken over all Viterbo cycles.
	\end{definition}
	
	\begin{proposition}\label{immersedtime}
		$\Sigma(t)$ cannot remain an embedding past time $T_{\text{imm}}=\frac{1}{4\pi}\Lambda^+$.
	\end{proposition}
	\begin{proof}
		This proposition follows from Viterbo's theorem and the fact that, for any Viterbo cycle $\gamma$,
			\begin{equation}
				\lambda(t).\gamma=\lambda(0).\gamma-4\pi t
			\end{equation}
	\end{proof}
	
	The similiarities between Proposition \ref{immersedtime} and Theorem C seem to indicate a connection between embeddedness and the type of the singularity. One can think of this connection as generalizing the case of plane curves, for which type II singularities only occur starting from nonembedded initial data. (In the case of an embedded plane curve, $\Lambda^+$ is twice the enclosed area, so for curves, Proposition \ref{immersedtime} is part of Grayson's Theorem.)
	
	\begin{conjecture}
		A Lagrangian mean curvature flow with embedded initial data which becomes immersed before the singular time must have a singularity of type II.
	\end{conjecture}
	
	\subsection{What are the minimal smooth blow-ups?}
Neves \cite{neveszeromaslov} has produced examples of Lagrangian surfaces in $\mathbf{C}^2$ whose singularities are of type II, and hence, have tangent flows whose supports consist of transversely-intersecting planes. His examples are given by rotating a curve $\gamma(s)\subset \mathbf{C}$: 
	\begin{equation}
		L_\gamma(r,s)=\left(\cos(r)\gamma(s),\sin(r)\gamma(s)\right)\in\mathbf{C}^2
	\end{equation}

It is a computation to see that the zero-Maslov condition requires the tangent flows given by Theorem \ref{TF.zeromaslov} to be $L_{\gamma_0}$, where $\gamma_0$ is the union of two lines intersecting at an angle of $\frac{\pi}{2}$. The smooth blow-up in Neves' examples is given by rotating an hyperbola asymptotic to these lines.

\begin{question}
	When $m=2$, is every minimal smooth blow-up of the form $L_\gamma$ for $\gamma$ an hyperbola?
\end{question}

\subsection{Are all type II singularities modeled by steady solitons?}
	The fact that type II smooth blow-ups are all eternal exact, zero Maslov solutions imposes a significant constraint on their geometry. The natural conjecture is that, just as type I smooth blow-ups are shrinking solitons, type II singularities ought to be steady solitons--namely, either minimal or translators.  But there has been no progress toward this conjecture to date. 

	In the case of curves in $\mathbf{C}$, Grayson showed \cite{graysonfigureeight} that certain type II singularities are modeled by the grim reaper translator. Joyce-Lee-Tsui \cite{joyce&lee&tsui10} have constructed translators which are asymptotic to planes; these translators are in fact exact. (It is an easy fact that all translators are automatically zero-Maslov.)

	\begin{question}
		Must a non-minimal smooth blow-up ust be a Joyce-Lee-Tsui translator, with asymptotics given by the blow-down cones?
	\end{question}

\bibliography{mcf}{}
\bibliographystyle{acm}

\end{document}